\documentclass[12pt]{article}

\usepackage{amsfonts,amsmath,amssymb}
\usepackage{mathrsfs,mathtools,url}
\usepackage{latexsym}
\usepackage{tikz,color}
\usepackage[T1]{fontenc}

\newtheorem{theorem}{Theorem}[section]

\newtheorem{lemma}[theorem]{Lemma}
\newtheorem{proposition}[theorem]{Proposition}

\newtheorem{corollary}[theorem]{Corollary}

\DeclareMathOperator {\diam} {diam}
\DeclareMathOperator {\rad} {rad}

\DeclareMathOperator {\Tr} {{\rm Tr}}

\def\cp{\,\square\,}
\newcommand{\proof}{\noindent{\bf Proof.\ }}
\newcommand{\qed}{\hfill $\square$ \bigskip}

\title{Extremal results on $k$-stepwise irregular graphs}

\author{Yaser Alizadeh$^{a}$, Sandi Klav\v{z}ar$^{b, c, d}$, Javaher Langari$^{a}$ \\\\
$^{a}$ \small Department of Mathematics, Hakim Sabzevari University, Sabzevar, Iran\\
\small {\tt y.alizadeh@hsu.ac.ir}\\
\small {\tt javaher\_langari61@yahoo.com}\\
$^{b}$ \small Faculty of Mathematics and Physics, University of Ljubljana, Slovenia\\
\small {\tt sandi.klavzar@fmf.uni-lj.si}\\
$^{c}$ \small Institute of Mathematics, Physics and Mechanics, Ljubljana, Slovenia \\
$^{d}$ \small Faculty of Natural Sciences and Mathematics, University of Maribor, Slovenia
}

\date{}
\begin{document}

\maketitle

\begin{abstract}
For a positive integer $k\ge 1$, a graph $G$ is $k$-stepwise irregular ($k$-SI graph) if the degrees of every pair of adjacent vertices differ by exactly $k$. Such graphs are necessarily bipartite. Using graph products it is demonstrated that for any $k\ge 1$ and any $d \ge 2$ there exists a $k$-SI graph of diameter $d$. A sharp upper bound for the maximum degree of a $k$-SI graph of a given order is proved. The size of $k$-SI graphs is bounded in general and in the special case when $\gcd(\Delta(G), k) = 1$. Along the way the degree complexity of a graph is introduced and used. 
\end{abstract}


\medskip\noindent
\textbf{AMS Math.\ Subj.\ Class.\ (2020)}: 05C07

%%%%%%%%%%%%%%%%%%%%%%%%%%%%%%%%%%%%%%%
\section{Introduction}
%%%%%%%%%%%%%%%%%%%%%%%%%%%%%%%%%%%%%%%

If $u$ and $v$ are two adjacent vertices of a graph $G = (V(G), E(G))$, then the {\em imbalance} of the edge $uv$ is  $|d_G(u) - d_G(v)|$, where $d_G(x)$ denotes the degree of the vertex $x$ of $G$. The {\em irregularity} of $G$ is then $\sum_{uv \in E(G)} |d_G(u)- d_G(v)|$~\cite{Albertson}. The study of graph irregularity has been a significant research direction. Many papers have been dedicated to investigating the measures of graph irregularity, a few examples are~\cite{Abdo1, Abdo2, Bickle}. The irregularity theory developed up to 2021 was then summarized in the book~\cite{Ali-2021}. Moreover, irregularity in graphs has been extensively utilized to analyze the topological structures and deterministic and random networks prevalent in chemistry, bio-informatics, and social networks, cf.~\cite{Estrada}.

A graph $G = (V(G), E(G))$ is a {\em stepwise irregular graph} (briefly {\em SI graph}) if for any edge $uv\in E(G)$ we have $|d_G(u) - d_G(v)| = 1$, where $d_G(x)$ denotes the degree of the vertex $x$ of $G$. This class of graphs was introduced by Gutman in~\cite{Gutman}, for further results on them, see~\cite{Adi-2022, Akhter-2022, Bera, Buya}. In particular, some sharp upper bounds on the maximum degree and maximum size of SI graphs were posed in~\cite{Bera, Buya}.  

A natural generalization of stepwise irregular graph form {\em $k$-stepwise irregular graphs} (briefly {\em $k$-SI graphs}) which are the graphs in which the two degrees of every pair of adjacent vertices differ by $k$, where $k$ is a fixed positive integer. This generalization was formally introduced in~\cite{Das-2023}, where it was mostly investigated for the case $k=2$. In this paper, we look in more detail at the general case and proceed as follows. In the next section, we first list some common definitions and concepts needed. Then we introduce the degree complexity of a graph and list a series of inequalities that hold in $k$-SI graphs. We also show that these graphs are bipartite. In Section~\ref{sec:max-degree} we use graph products to demonstrate that for any $k\ge 1$ and any $d \ge 2$ there exists a $k$-SI graph of diameter $d$. In the subsequent section we give a sharp upper bound for the maximum degree of a $k$-SI graph of a given order. In Section~\ref{sec:size} we bound the size of $k$-SI graphs in general and in the special case when $\gcd(\Delta(G), k) = 1$.

%%%%%%%%%%%%%%%%%%%%%%
\section{Preliminaries}
%%%%%%%%%%%%%%%%%%%%%%

Let $G = (V(G), E(G))$ be a graph. The order and the size of $G$ will be denoted by $n(G)$ and $m(G)$. The degree of $v\in V(G)$ is denoted by $d_G(v)$, and $\Delta(G)$ and $\delta(G)$ denote the maximum and the minimum degree in $G$, respectively. The {\em distance} $d_G(u,v)$ between vertices $u$ and $v$ of a graph $G$ is the minimum length between the $u,v$-paths. The {\em eccentricity} $\varepsilon_G(v)$ of a vertex $v$ is the maximum distance between $v$ and the other vertices in $G$. The {\em diameter} $\diam(G)$ and the {\em radius} $\rad(G)$ are the maximum and the minimum eccentricity of vertices in $G$, respectively. For a positive integer $k$ we will use the notation $[k] = \{1,\ldots, k\}$. 

We introduce the {\em degree complexity} $C_d(G)$ of $G$ as the number of distinct degrees in $G$. (For a general approach to the concept of the complexity of a graph invariant, see for example~\cite{alizadeh-2016}.) Note that $C_d(G) = 1$ if and only if $G$ is regular. Note also that every graph has at least two vertices of the same degree, hence $C_d(G) \le n(G) - 1$. 

For a $k$-SI graph $G$ we set 
$$A_i = \{u\in V(G):\ d_G(u) = \Delta(G) - ik\},\ i\in \{0, 1, \ldots, C_d(G) -1\}\,.$$
Since $G$ is a $k$-SI graph, the sets $A_i$ are nonempty, pairwise disjoint, and each of them induces an edgeless graph. We now derive inequalities on $a_r = |A_r|$, which will be useful in the rest of the paper.  

Since the sets $A_i$ are pairwise disjoint, $n(G) = \sum_{r=0}^{C_d(G)-1} a_r$, where $C_d(G) = \frac{\Delta - \delta}{k} + 1$. Further, if $u \in A_0$ and $v \in A_{C_d(G)-1}$, then a vertex adjacent to $u$ is of degree $\Delta(G)-k$, and a vertex adjacent to $v$ is of degree $\delta(G) + k$. That is, $N_G(u)\subseteq A_1$ and $N_G(v)\subseteq A_{C_d(G)-2}$, which in turn implies 
\begin{equation}\label{eq1}
    a_1\geq \Delta(G)
\end{equation}
and 
\begin{equation}\label{eq2}
a_{C_d(G)-2}\geq \delta(G)\,.
\end{equation}
Similarly, if $i \in [C_d(G)-2]$ and $w \in A_i$, then $ N(w) \subseteq A_{i-1} \cup A_{i+1}$, which yields 
\begin{equation}\label{eq3}
\Delta(G) -ik \leq a_{i-1}+a_{i+1},\ i \in [C_d(G)-2]\,.
\end{equation}

For $0 \leq i \leq C_d(G)-1$, let $E(A_i)$ denote the set of edges incident with a vertex in $A_i$. Since $A_i$ is an independent set, we have $|E(A_i)| = a_i (\Delta-ik)$. Furthermore, by definition, $E(A_i ) \subseteq E(A_{i-1} \cup A_{i +1} )$ for $i\in [C_d(G)-1]$, hence 
\begin{equation}\label{eq4}
   a_i( \Delta(G) -ik) \leq a_{i-1}(\Delta(G) - (i-1)k) + a_{i+1}(\Delta(G) - (i+1)k)\,. 
\end{equation}
For the vertices of the minimum and the maximum degree we respectively have $E(A_0) \subseteq E(A_1)$ and $E(A_{C_d(G)-1}) \subseteq E(A_{C_d(G)-2})$, therefore 
\begin{equation}\label{eq5}
  a_0 \Delta(G)  \leq a_1 (\Delta(G)-k)
\end{equation}
and 
\begin{equation}\label{eq6}
  a_{C_d(G)-1} \delta(G) \leq a_{C_d(G)-2} (\delta(G)+k)\,.
\end{equation}
Note that the equalities in~\eqref{eq5} and~\eqref{eq6} hold if and only if $C_d(G)=2$.

By~\eqref{eq3} we have $a_0 + a_2 \ge \Delta(G)-k$. If $a_0 + a_2=\Delta(G)-k$, then 
\begin{align*}
\Delta(G)(\Delta(G)-k) & \leq a_1 (\Delta(G)-k) \qquad\qquad\qquad\quad ({\rm by}~\eqref{eq1})\\
& \leq a_0 \Delta(G) + a_2 (\Delta(G)-2k) \qquad ({\rm by}~\eqref{eq4}) \\ 
& = \Delta(G)(a_0+a_2)-2ka_2 \\
& = \Delta(G)(\Delta(G)-k)-2a_2k\,.
\end{align*}
We can conclude that if $a_0 + a_2=\Delta(G)-k$, then $a_2=0$. Hence, if $a_2>0$ (equivalently $C_d(G)\ge 3$), then 
\begin{equation}\label{eq7}
  a_0 + a_2 \geq \Delta(G)-k + 1\,.
\end{equation}

To conclude the preliminaries, we prove the next basic result on $k$-SI graphs which was for the case $k=1$ established in~\cite[Lemma~3]{Gutman}, and for the case $k=2$ in~\cite[Theorem 3]{Das-2023}. 

\begin{proposition}
\label{prop:bipartite} 
If $k\ge 1$ and $G$ is a $k$-SI graph, then $G$ is bipartite. 
\end{proposition}

\proof 
Suppose on the contrary that $G$ contains on odd cycle whose consecutive vertices are $v_0, v_1,\ldots, v_{2t}$. We may assume that $d_G(v_0)=\min \{d_G(v_i):\ 0 \le i \le 2t \}$. Since $G$ is a $k$-SI graph, $d_G(v_i) \equiv d_G(v_{2t-i+1}) \bmod{2k}$. Thus $v_t \equiv v_{t+1} \bmod{2k}$, but this is not possible as  $v_tv_{t+1}\in E(G)$. 
\qed

\begin{corollary}
The unique $k$-SI graph $G$ of diameter $2$ is $G\cong K_{\frac{n(G)+k}{2},\frac{n(G)-k}{2}}$. 
\end{corollary}

\proof 
Let $G$ be a $k$-SI graph with $\diam(G) = 2$. By Proposition~\ref{prop:bipartite}, the graph $G$ is bipartite, which in turn implies that $G$ is a complete bipartite graph. As $G$ is a $k$-SI graph, this means that $G \cong K_{m, m+k}$ for some integer $m$. Therefore, This implies that $n(G) = 2m(G) + k$, so that $G \cong K_{\frac{n+k}{2},\frac{n-k}{2}}$. 
\qed

%%%%%%%%%%%%%%%%%%%%%%%%%%%%%
\section{Constructing $k$-SI graphs using graph produces}
\label{sec:constructions}
%%%%%%%%%%%%%%%%%%%%%%%%%%%%%

When dealing with new mathematical structures, it is a desirable to first demonstrate that they actually exist. Several examples of $k$-SI graphs have been already presented in~\cite{Das-2023}, with a focus on 2-SI graphs. In this section we show that for any $k\ge 2$ and any diameter $d\ge 2$ there are graphs $G$ with $\diam(G) = d$ which are  $k$-SI graphs. For this sake we will use Cartesian and lexicographic products of graphs. Before we do that, let's look at two sporadic families of such graphs which are introduced in~\cite{ali-2024} for a different purpose. 

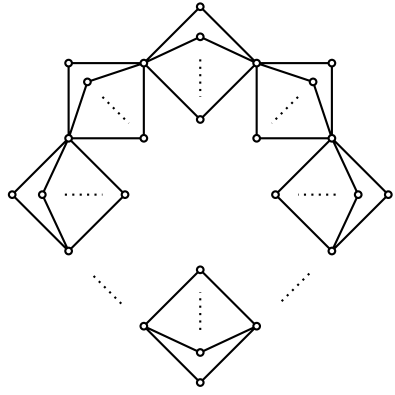
\begin{figure}[ht!]
\begin{center}
\begin{tikzpicture}[scale=0.5,style=thick]
\tikzstyle{every node}=[draw=none,fill=none]
\def\vr{2.5pt} 

\begin{scope}[yshift = 0cm, xshift = 0cm]
\path (7.5,9.5) coordinate (v0);
\path (9.5,7.5) coordinate (v2);
\path (9.5,4.5) coordinate (v3);
\path (7.5,2.5) coordinate (v4);
\path (4.5,2.5) coordinate (v5);
%%v6&\ddots&v5
\path (2.5,4.5) coordinate (v6);
\path (2.5,7.5) coordinate (v7);
\path (4.5,9.5) coordinate (v1);
\path (9.5,9.5) coordinate (u1);
\path (11,6) coordinate (u2);
\path (6,1) coordinate (u3);
\path (1,6) coordinate (u4);
\path (2.5,9.5) coordinate (u5);
\path (6,11) coordinate (w{1,1});
\path (7.5,7.5) coordinate (w1);
\path (8,6) coordinate (w2);
\path (6,4) coordinate (w3);
\path (4,6) coordinate (w4);
\path (4.5,7.5) coordinate (w5);
\path (6,8) coordinate (w{1,p});

\path (9,9) coordinate (z1);
\path (10.2,6) coordinate (z2);
\path (6,1.8) coordinate (z3);
\path (1.8,6) coordinate (z4);
\path (3,9) coordinate (z5);
\path (6,10.2) coordinate (w{1,2});
%% edges %%
\draw (v6) -- (z4) -- (v7) -- (z5) -- (v1) -- (w{1,2}) -- (v0) -- (u1) -- (v2) -- (u2) -- (v3) -- (w2) -- (v2) -- (w1) -- (v0) --(w{1,p}) -- (v1) -- (w5) -- (v7) -- (w4) -- (v6) -- (u4) -- (v7) -- (u5) -- (v1) -- (w{1,1}) -- (v0) -- (z1) -- (v2) --(z2) -- (v3);  
\draw (v4) -- (u3) -- (v5) -- (w3) -- (v4) -- (z3) -- (v5);
%%\draw  (8.5, 3.5) circle (1cm);
\draw [dotted] (8.6,8.6) -- (7.9,7.9);
\draw [dotted] (9.6,6) -- (8.6,6);
\draw [dotted] (6,2.4) -- (6,3.4);
\draw [dotted] (2.4,6) -- (3.4,6);
\draw [dotted] (3.4,8.6) -- (4.1,7.9);
\draw [dotted] (6,9.6) -- (6,8.6);
\draw [dotted] (8.9,3.9) -- (8.1,3.1);
\draw [dotted] (3.9,3.1) -- (3.1,3.9);
%% vertices %%%

%%\draw (1,2) .. controls (.2.5,3) .. (1,4);
\draw (v0)  [fill=white] circle (\vr);
\draw (v2)  [fill=white] circle (\vr);
\draw (v3)  [fill=white] circle (\vr);
\draw (v4)  [fill=white] circle (\vr);
\draw (v5)  [fill=white] circle (\vr);
\draw (v6)  [fill=white] circle (\vr);
\draw (v7)  [fill=white] circle (\vr);
\draw (v1)  [fill=white] circle (\vr);
\draw (u1)  [fill=white] circle (\vr);
\draw (u2)  [fill=white] circle (\vr);
\draw (u3)  [fill=white] circle (\vr);
\draw (u4)  [fill=white] circle (\vr);
\draw (u5)  [fill=white] circle (\vr);
\draw (w{1,1})  [fill=white] circle (\vr);
\draw (w1)  [fill=white] circle (\vr);
\draw (w2)  [fill=white] circle (\vr);
\draw (w3)  [fill=white] circle (\vr);
\draw (w4)  [fill=white] circle (\vr);
\draw (w5)  [fill=white] circle (\vr);
\draw (w{1,p})  [fill=white] circle (\vr);

\draw (z1)  [fill=white] circle (\vr);
\draw (z2)  [fill=white] circle (\vr);
\draw (z3)  [fill=white] circle (\vr);
\draw (z4)  [fill=white] circle (\vr);
\draw (z5)  [fill=white] circle (\vr);
\draw (w{1,2})  [fill=white] circle (\vr);

\end{scope}
\end{tikzpicture}
\end{center}
\caption{The graph $\Gamma_{p,q}$ }
\label{Gpq}
\end{figure}

\begin{figure}[ht!]
\begin{center}
\begin{tikzpicture}[scale=0.7,style=thick]
\tikzstyle{every node}=[draw=none,fill=none]
\def\vr{2.5pt} 

\begin{scope}[yshift = 0cm, xshift = 0cm]
\path (7.5,9.5) coordinate (v1);
\path (9.5,8) coordinate (v2);
\path (11,6) coordinate (v3);
\path (9.5,4.8) coordinate (v4);
\path (7.5,3.3) coordinate (v5);
\path (5.5,4.8) coordinate (v6);
\path (4,6) coordinate (v7);
\path (5.5,8) coordinate (v8);
\path (9.3,9.5) coordinate (u1);
\path (11,7.8) coordinate (u2);
\path (9.3,3.3) coordinate (u3);
\path (5.7,3.3) coordinate (u4);
\path (4,7.8) coordinate (u5);
\path (5.7,9.5) coordinate (u6);
\path (7.9,8) coordinate (w1);
\path (9.5,6.4) coordinate (w2);
\path (7.9,4.8) coordinate (w3);
\path (7.1,4.8) coordinate (w4);
\path (5.3,6.3) coordinate (w5);
\path (7.1,8) coordinate (w6);

\path (8.9,9.1) coordinate (z1);
\path (8.9,3.7) coordinate (z2);
\path (4.4,7.3) coordinate (z3);

%% edges %%
\draw (v1) -- (u1) -- (v2) -- (u2) -- (v3) -- (w2) -- (v2) -- (w1) -- (v1) -- (w6) -- (v8) -- (w5) -- (v7) -- (u5) -- (v8) --(u6) -- (v1); %%   
\draw (v4) -- (u3) -- (v5) -- (u4) -- (v6) -- (w4) -- (v5) -- (w3) -- (v4);
\draw (v1) -- (z1) -- (v2);
\draw (v4) -- (z2) -- (v5);
\draw (v7) -- (z3) -- (v8);

\draw  [dotted](10.6,5.5) -- (9.9,5.1);
\draw  [dotted](5.1,5.1) -- (4.4,5.7);
\draw  [dotted](8.6,8.8) -- (8.2,8.3);
\draw  [dotted](8.6,4.0) -- (8.1,4.5);
\draw  [dotted](4.6,7.05) -- (5.1,6.55);

%% vertices %%%

\draw (v1)  [fill=white] circle (\vr);
\draw (v2)  [fill=white] circle (\vr);
\draw (v3)  [fill=white] circle (\vr);
\draw (v4)  [fill=white] circle (\vr);
\draw (v5)  [fill=white] circle (\vr);
\draw (v6)  [fill=white] circle (\vr);
\draw (v7)  [fill=white] circle (\vr);
\draw (v8)  [fill=white] circle (\vr);
\draw (u1)  [fill=white] circle (\vr);
\draw (u2)  [fill=white] circle (\vr);
\draw (u3)  [fill=white] circle (\vr);
\draw (u4)  [fill=white] circle (\vr);
\draw (u5)  [fill=white] circle (\vr);
\draw (u6)  [fill=white] circle (\vr);
\draw (w1)  [fill=white] circle (\vr);
\draw (w2)  [fill=white] circle (\vr);
\draw (w3)  [fill=white] circle (\vr);
\draw (w4)  [fill=white] circle (\vr);
\draw (w5)  [fill=white] circle (\vr);
\draw (w6)  [fill=white] circle (\vr);

\draw (z1)  [fill=white] circle (\vr);
\draw (z2)  [fill=white] circle (\vr);
\draw (z3)  [fill=white] circle (\vr);

\end{scope}

\end{tikzpicture}
\end{center}
\caption{The graph $H_{p,q}$ }
\label{Hpq}
\end{figure}
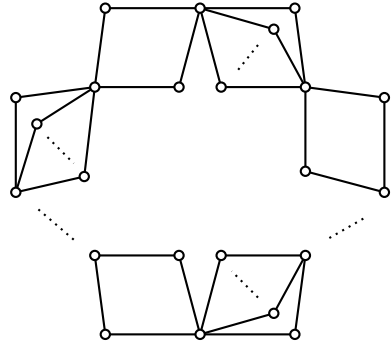

Let $p,q \ge 2$. Then the graph $\Gamma_{p,q}$ consists of $q$ induced subgraphs $K_{2,p}$ arranged in a cyclic structure as shown in Fig.~\ref{Gpq}. The graph $H_{p,q}$ consists of $q$ induced subgraphs $K_{2,p}$ and $q$ induced subgraphs $K_{2,2}$ alternately arranged in a cyclic structure as shown in Fig.~\ref{Hpq}. Then $\Gamma_{p,q}$ is a $(2p-2)$-SI graph of diameter $q$ and $H_{p,q}$ is a $p$-SI graph of diameter $2q$. 

%%%%%%%%%%%%%%%%%%%%%%%%%%%%%%%%%%%%%%%%

Recall that the {\em Cartesian product} $G\cp H$ of two graphs $G$ and $H$ is the graph with $V(G\cp H) = V(G)\times V(H)$, along with the condition that two vertices $(g, h)$ and $(g', h')\in V(G) \times V(H)$ are adjacent if $g = g'$ and $hh' \in E(H)$, or $gg' \in E(G)$ and $h = h'$. In~\cite[Theorem~5]{Das-2023} it was shown that the Cartesian product of two $2$-SI graphs is a $2$-SI graph; moreover, it was left as an exercise to the reader that the same holds for two $k$-SI graphs. We show that a more general statement holds. 

\begin{proposition}
\label{prop:cartesian}
If $G$ and $H$ are graphs, then $G\cp H$ is a $k$-SI graph if and only if both $G$ and $H$ are $k$-SI graphs.
\end{proposition}

\proof  
Assume that $G\cp H$ is a $k$-SI graph. Let $gg' \in E(G)$ and $h \in V(H)$. Since $G\cp H$ is a $k$-SI graph and by the fact $d_{G\cp H}((g,h))= d_G(g) + d_H(h)$, we have
\[|d_{G\cp H}((g,h)) - d_{G\cp H}((g',h))|= |d_G(g) -d_G(g')| = k.\]
Thus $G$ is a $k$-SI graph. Analogously, $H$ is a $k$-SI graph.

Assume now that $G$ and $H$ are $k$-SI graphs and let $(g,h)(g,h') \in E(G\cp H)$. Since $H$ is a $k$-SI graph, we have $k = |d_H(h) - d_H(h')| = |d_{G\cp H}((g,h)) - d_{G\cp H}((g,h'))|$. Similarly, since $G$ is a $k$-SI graph, for an edge $(g,h)(g',h) \in E(G\cp H)$ we have 
$k = |d_G(g) - d_G(g')| = |d_{G\cp H}((g,h)) - d_{G\cp H}((g',h))|$.
\qed

We next consider the {\em lexicographic product} $G\circ H$ of two graphs $G$ and $H$ which also has $V(G\cp H) = V(G)\times V(H)$, and  $(g, h)(g', h')\in E(G) \circ V(H)$ if $g = g'$ and $hh' \in E(H)$, or $gg' \in E(G)$. We recall that if $G$ is not complete (and connected), then $\diam(G\circ H) = \diam(G)$~\cite{HIK-2011}. Denoting by $\overline{K}_t$ the complement of $K_t$, we have the following result. 

\begin{proposition}
\label{prop:lex}
If $G$ is a $k$-SI graph, then $G\circ \overline{K}_t$ is a $(kt)$-SI graph. 
\end{proposition}

\proof  
Since the second factor of the lexicographic product considered is edgeless, an arbitrary edge of $G\circ \overline{K}_t$ is of the form $(g,h)(g',h')$, where $gg'\in E(G)$ and $h,h'\in V(\overline{K}_t)$. Then, by the definition of the lexicographic product and the assumption that $G$ is a $k$-SI graph, we have 
\begin{align*}
|d_{G\circ H}((g,h)) - d_{G\cp H}((g',h'))| & =  
|t\cdot d_G(g) - t\cdot d_G(g')| \\ 
& = t\cdot |d_G(g) - t\cdot d_G(g')| = kt\,,    
\end{align*}
which is exactly what we wanted to see.
\qed

\begin{theorem}
\label{thm:all-cases-possible}
For any $k\ge 1$ and any $d \ge 2$, there exists a $k$-SI graph of diameter $d$.
\end{theorem}

\proof 
If $m\ge 2$, then the complete bipartite graph $K_{m,m+k}$ is a $k$-SI graph for any $k\ge 1$. This settles the case $d=2$ and all $k\ge 1$.  

Let next $G$ be a $1$-ST graph of diameter $3$, say the graph presented in~\cite[Fig.~1]{Buya}. Then by Proposition~\ref{prop:lex}, the lexicographic product $G\circ \overline{K}_k$ is a $k$-SI graph for any $t\ge 1$. Since $\diam(G\circ \overline{K}_k) = 3$, this settles the case $d=3$ and all $k\ge 1$. 

Proceeding by induction, let $k\ge 1$ and $d\ge 4$. Select an arbitrary $k$-SI graph $G$ with $\diam(G) = 2$ and an arbitrary $k$-SI graph $H$ with $\diam(H) = d-2$. Then by Proposition~\ref{prop:cartesian}, the Cartesian product $G\cp H$ is a $k$-SI graph. Since $\diam(G\cp H) = \diam(G) + \diam(H) = d$ we are done. 
\qed

Note that the proof method of Theorem~\ref{thm:all-cases-possible} can be used to obtain different infinite families of $k$-SI graphs. Just as an example, if $k\ge 1$ and $t\ge 2$ are fixed, then $K_{m,m+k}\circ \overline{K}_t$ is a $(kt)$-SI graph of diameter $2$ for any $m\ge 2$. 

However, for certain classes of graphs, not all diameters are realizable, as the following result shows. 

%%%%%%%%%%%%%%%%%%%%%%%%%%%%%
\begin{proposition}
\label{even-diam}
$k$-SI trees and $k$-SI unicyclic graphs are of even diameter.
\end{proposition}

\proof
Suppose that $T$ is a $k$-SI tree of odd diameter with diametrical path on the vertices $v_1, v_2, \ldots v_{2t}$. Clearly, $d_T(v_1) = d_T(v_{2t}) = 1$. Since $T$ is a $k$-SI graph, $d(v_i) \equiv d(v_{2t-i+1}) \bmod 2k$. This implies that $d_T(v_t) \equiv d_T(v_{t+1}) \bmod 2k$, a contradiction.  

Let next $G$ be a $k$-SI unicyclic graph, and let $C$ be its unique cycle $C$. Let $P$ be a diametrical path of $G$. If $C$ contains none of the end vertices of $P$, then end vertices of $P$ are of degree $1$. By the above argument we get that the length of $P$ is even.  Otherwise, one of the end vertices of $P$, say $u$, is located on $C$. If $d_G(u) = 2$, then the two neighbours of $u$ are of degree $k+2$. Thus we can find another diametrical path whose both end vertices are of degree $1$ and we can conclude as above that $\diam(G)$ is even. Finally, if $d_T(u) \ge 3$, then we can prolong $P$ with a vertex not on $C$ and adjacent to $u$. So this case cannot happen. 
\qed

%%%%%%%%%%%%%%%%%%%%%%%%%%%%%%%%%%%%%%%%%%%%%%%%%%%%%%%%%%
\section{Bounding maximum degree in $k$-SI graphs}
\label{sec:max-degree}
%%%%%%%%%%%%%%%%%%%%%%%%%%%%%%%%%%%%%%%%%%%%%%%%%%%%%%%%%%

In this section we give a sharp upper bound on the maximum degree of a $k$-SI graph of a given order. To prove it, we will use the findings from the proof of the next result. For it we recall that if $G$ is a $k$-SI graph, then $A_i$ is the set of vertices of $G$ of degree $\Delta(G) - ik$ and that $a_i = |A_i|$. 

\begin{proposition}
Let $G$ be a $k$-SI graph. If $k\equiv \Delta(G) \bmod 2$, then $m(G)$ is even.
\end{proposition}

\proof
By Proposition~\ref{prop:bipartite}, $G$ is bipartite. Let $X, Y$ be the bipartition of $G$, where $X$ contains a maximum degree vertex. If $\Delta(G)$ is even, then the evenness of $m(G)$ follows from the equality 
\begin{align}
\label{eq8}
m(G) = & \underset{{u\in X}} \sum d_G(u)=\underset{{i\geq 0}}\sum a_{2i}(\Delta(G)-(2i)k)\,. 
\end{align}
On the other hand, we also have the equality
\begin{align}
\label{eq9}
m(G) = & \underset{{v\in Y}}\sum d_v(v)=\underset{{i\geq 0}}\sum a_{2i+1}(\Delta(G)-(2i+1)k)\,,
\end{align}
from which the evenness of $m(G)$ follows when $\Delta(G)$ is odd. 
\qed

The main result of this section bounds the maximum degree of $k$-SI graphs and reads as follows. 

\begin{theorem}
\label{max degree}
If $G$ is a $k$-SI graph, then  
$$\Delta(G) \leq \left\lfloor \frac{n(G)+k}{2} \right\rfloor\,,$$
where the equality holds if and only if $G \cong K_{\frac{n(G)+k}{2}, \frac{n(G)-k}{2}}$. 
\end{theorem}

\proof 
Using~\eqref{eq2} we have
\begin{align}
m(G) = &\  \sum _{i\geq 0}a_{2i+1}(\Delta(G)-(2i+1)k)\nonumber \\
= &\  a_1 (\Delta(G) - k)+ \sum _{i\geq 1}a_{2i+1}(\Delta(G)-(2i+1)k) \nonumber \\
\geq &\  a_1(\Delta(G)-k). \label{eq11}
\end{align}  
Using~\eqref{eq1} and~\eqref{eq2} and the fact $\sum _{i\geq 0}^{C_d(G)-1} a_{i}=n(G)$, we further get
\begin{align}
2m(G) = &\ \sum _{i\geq 0}^{C_d(G)-1} a_{i}(\Delta(G)-ik) = \Delta(G) \sum_{i\ge 0}a_i - ka_1-  \sum _{i= 2}^{C_d(G)-1} ika_i \nonumber \\ 
& \leq \Delta(G) n(G) - ka_1. \label{eq12}
\end{align} 
Hence, from~\eqref{eq11}, \eqref{eq12}, and~\eqref{eq1} we obtain
\begin{equation}\label{eq13}
n(G)\Delta(G) \geq a_1(2\Delta(G)-k)\geq \Delta(G)(2\Delta(G)-k)\,,
\end{equation} 
which implies $\Delta(G) \leq \left\lfloor \frac{n(G)+k}{2} \right\rfloor$. This proves the theorem's inequality. Assume now that the equality holds. We distinguish two cases.

\medskip\noindent
\textbf{Case 1}: $n(G)\equiv k \bmod 2$. \\
In this case, $\Delta(G) = \frac{n(G)+k}{2}$. Since the equalities in~\eqref{eq11},~\eqref{eq12} and~\eqref{eq13} must hold, we infer that $a_1=\Delta(G)$, $a_2 = 0$, and $a_0 = n(G) - \Delta(G) = \frac{n(G)-k}{2}$. On the other hand, since the equality in~\eqref{eq11} holds and $\Delta(G) = \frac{n(G)+k}{2}$, we get $G\cong K_{\frac{n(G)+k}{2}, \frac{n(G)-k}{2}}$.

\medskip\noindent
\textbf{Case 2}: $n(G)$ and $k$ have different parities. \\
Now $\Delta(G) = \frac{n(G)+k-1}{2}$. If $a_1 \ge \Delta(G) + 1$, then by~\eqref{eq13},  
\begin{align*}
(2\Delta(G) -k + 1) \Delta(G) = & \ n(G)\Delta(G) \ge a_1(2\Delta(G) - k)\\
& \ge \ (\Delta(G) + 1)(2\Delta(G)-k+1)\,,
\end{align*}
and hence $\Delta(G) \le k$, a contradiction. Thus $a_1 = \Delta(G)$. If $a_2 = 0$, then since $a_0 + a_1 = n(G)$, we get $a_0 = \Delta(G) - k + 1$. By~\eqref{eq8} and~\eqref{eq9} it follows that 
\[(\Delta(G) - k +1) \Delta(G) = a_0 \Delta(G) = m(G)= a_1(\Delta(G) - k) = \Delta(G)(\Delta(G) - k)\,, \] 
another contradiction. Therefore $a_2 >0$ which in turn implies that $G$ is not a complete bipartite graph, so we are done.
\qed

%%%%%%%%%%%%%%%%%%%%%%%%%%%%%%%%%%%%%%%%%%%%%%%%%%%%%%%%%%
\section{Bounding size in $k$-SI graphs}
\label{sec:size}
%%%%%%%%%%%%%%%%%%%%%%%%%%%%%%%%%%%%%%%%%%%%%%%%%%%%%%%%%%

In this section we bound the size of $k$-SI graphs. In our first main result we give a sharp general upper bound, and in the second main result we give an upper bound for the case when $\gcd(\Delta(G), k) = 1$.

\begin{theorem}\label{max:size1}
If $G$ is a $k$-SI graph, then
$$m(G )\leq \frac{n\Delta(G)(\Delta(G)-k)}{2\Delta-k}\,,$$
where the equality holds if and only if $C_d(G)=2$.
\end{theorem}

\proof 
Combining~\eqref{eq8} and~\eqref{eq9} we get
\[2m(G) = n(G)(\Delta(G) - k) + a_0k -\Big(ka_2 + 2ka_3 + \cdots +(t-1)ka_{t}\Big)\,,\]
where $t=C_d(G)-1$. It follows that $2 m(G) \leq n(G) (\Delta(G) - k) + a_0k$. Hence by~\eqref{eq5},      
\begin{equation}\label{eq22}
    2 m(G) \leq (\Delta(G) -k)\left(n(G) + \frac{a_1k}{\Delta(G)}\right)\,.
\end{equation}
Using~\eqref{eq4} we get $a_0 + a_1 + a_2 \leq n(G)$ and $\Delta(G) a_0 +(\Delta(G) -2k)a_2 \geq (\Delta(G) -k )a_1$. Using these inequalities we can estimate as follows: 
\begin{align*}
\Delta(G)(n(G)-a_1) \geq & \ \Delta(G)(a_0+a_2) \geq \Delta(G)(a_0 +a_2)-2k a_2 \\
= & \ \Delta(G) a_0 +(\Delta(G)-2k)a_2 \\
\geq & \  (\Delta(G) -k )a_1\,.
\end{align*}
From here we get
\begin{equation}\label{eq23}
    a_1\leq \frac{n(G)\Delta(G)}{2\Delta(G)-k}\,.
\end{equation}
Combining~\eqref{eq22} with~\eqref{eq23} we get the required inequality. 

The equalities in~\eqref{eq22} and in~\eqref{eq23} hold if and only if $a_2 = 0$, that is, if and only if $C_d(G) =2$. 
\qed

In view of Proposition~\ref{prop:bipartite}, the only $k$-SI graphs of radius $1$ are stars. In the next consequence of Theorem~\ref{max:size1} we may thus restrict to graphs of radius at least $2$. To state the corollary, we recall a couple of well-established concepts.

The {\em transmission} $\Tr_G(v)$ of a vertex $v\in V(G)$  is the sum of distances between $v$ and the other vertices of $G$. The transmission of a vertex is one of the basic concepts in metric graph theory, cf.~\cite{plenski, soltes}. We note in passing that stepwise transmission irregular graphs were extensively considered~\cite{ali-2024, al-yakoob-2022, damn-2024, damn-2024b, dobrynin-2020, xu-2023}. The celebrated Wiener index~\cite{wiener} $W(G)$ of $G$ can then be expressed as $W(G) = \frac{1}{2}\sum _{v\in V(G)}\Tr_G(v)$. Finally, $G$ is {\em $d$-self-centered} if $\varepsilon_G(v)$ holds for all vertices of $G$.

\begin{corollary}\label{wiener}
If $G$ is a $k$-SI graph with $\rad(G) \ge 2$, then
\[ W(G) \ge n(G)\left[ n(G) - \frac{\Delta(G)(\Delta(G)-k)}{2\Delta(G)-k} - 1 \right]\,.\]
Moreover, the equality holds if and only if $G$ is a $2$-self centered graph with $C_d(G) = 2$. 
\end{corollary}

\proof
If $v\in V(G)$, then since $\varepsilon(v) \ge 2$, we have 
\begin{equation}
\Tr_G(v) \ge d_G(v) + 2(n(G)-1-d_G(v))= 2(n(G)-1) - d_G(v)\,.
\label{eq:2-self}
\end{equation}
Summing over all vertices and having in mind that $2W(G) = \sum _{v\in V(G)}\Tr_G(v)$ we get
\[W(G) \ge n(G)(n(G)-1) - m(G)\,.\]
The desired inequality now follows from Theorem~\ref{max:size1}. 

The equality in~\eqref{eq:2-self} holds if and only if $G$ is a $2$-self centered graph, while the equality in Theorem~\ref{max:size1} holds if and only if $C_d(G) = 2$, from which the equality assertion of the corollary follows. 
\qed

If $k\ge 2$ and $n\ge 2$, then $K_{n+k,n}$ is a $k$-SI graph which fulfills the conditions for equality in Corollary~\ref{wiener}.   

To prove the announced second main result, we need the following auxiliary result. 

\begin{lemma}
\label{prime}
If $G$ is a $k$-SI graph, $C_d(G)=2$, and $\gcd(\Delta(G), k) =1$, then $2\Delta(G) -k$ divides $n(G)$.  
\end{lemma}

\proof
$C_d(G)=2$ implies $a_2 =0$. Thus $a_0 + a_1 =n(G)$ and $\Delta(G) a_0 = (\Delta(G) -k)a_1$. It follows that 
\[a_1 = \frac{n\Delta(G)}{2\Delta(G) -k}\,.\]
Since $\Delta(G)$ and $2\Delta(G) -k$ are co-prime, $2\Delta(G)-k$ divides $n(G)$.
\qed

\begin{theorem}
\label{max:size2}
If $G$ is a $k$-SI graph,  $G\ne K_{\Delta(G), \Delta(G)-k}$, and $\gcd(\Delta(G), k) = 1$, then 
\begin{equation}\label{eq24}
    m(G) \leq (\Delta(G) - k)(n(G)-\Delta(G) +k-1)\,.
\end{equation}
Moreover, equality holds if and only if $C_d(G)=3$ and $|X| = \Delta(G)-k+1$, where $X$ is the bipartition part of $G$ that contains a maximum degree vertex.
\end{theorem}

\proof  
Assume first that $C_d(G)=2$. Then $a_2 = 0$ and $2\Delta(G)-k$ divides $n(G)$
by Lemma~\ref{prime}. It follows that $2\Delta(G) - k\leq n(G)$. If $n(G) = 2\Delta(G)-k$, then Theorem~\ref{max degree} implies that $G$ is isomorphic to $K_{\Delta(G),\Delta(G)-k}$, which we have excluded in the theorem's statement.  Hence $2\Delta(G)-k \leq \frac{n(G)}{2}$, and applying Theorem~\ref{max:size1} we obtain
\begin{equation}\label{eq27}
\begin{array}{rllll}
m(G) &=& \dfrac{n(G)\Delta(G)(\Delta(G)-k)}{2\Delta(G)-k} = \left( n(G) - \dfrac{n(G)(\Delta(G) -k)}{2\Delta(G)-k} \right)(\Delta(G) -k)\vspace{.3cm}\\
&\leq & \left( n(G) - \dfrac{(4\Delta(G) - 2k)(\Delta(G) -k)}{2\Delta(G) -k} \right)(\Delta(G)-k) \vspace{.3cm}\\
&= &(n(G) - 2\Delta(G) + 2k)(\Delta(G) - k) \vspace{.3cm}\\
&\leq & (n(G) - \Delta(G) + k - 1)(\Delta(G) - k)\,,
\end{array}
\end{equation}
where we have used that $\Delta(G) \geq k+1$. If the equalities in~\eqref{eq27} hold, then we have $n(G) = 4\Delta(G) -2k$ and $\Delta(G)=k+1$. Hence $\delta(G) = \Delta(G) -k = 1$. Therefore, in this case, $G$ is isomorphic to  $K_{\Delta,1}$. As this was excluded, we can conclude that $C_d(G)=2$ is not possible. 

In the rest we may thus assume that $C_d(G)\ge 3$. Hence $a_2 > 0$ and $a_0 + a_2 \geq \Delta(G)-k+1$ by~\eqref{eq7}. Then
\begin{align}
m(G) \geq & \ a_0 \Delta(G) + a_2 (\Delta(G) - 2k) = (a_0 + a_2 )(\Delta(G) - k) + a_0k -a_2k \nonumber \\ 
\geq & \ (\Delta(G)-k)(\Delta(G)-k+1) + a_0k - a_2k\,. \label{eq25}
\end{align}
Moreover, we have
\begin{align*}
2 m(G) = & \ a_0\Delta(G)+ a_1 (\Delta(G)-k) +\cdots + a_{t}\delta(G) \\
& = \ n (\Delta(G)-k) + a_0k -a_2k -2 a_3k - \cdots -(t-1)a_tk\,,
\end{align*}
where $ t = C_d(G)-1$. Thus
\begin{equation}\label{eq26}
    2 m(G) \leq n(G) (\Delta(G)-k) + a_0k - a_2k\,.
\end{equation}
From~\eqref{eq25} and~\eqref{eq26} we get the required inequality. The equalities in~\eqref{eq25} and~\eqref{eq26} hold if and only if $a_3 = 0$ and $a_0 + a_2 = \Delta(G) -k+1$. Thus, the equality in~\eqref{eq24} holds if and only if $C_d(G)=3$ and $|X| = a_0 + a_2= \Delta(G) - k+1$ because $a_2 > 0$.
\qed

\section*{Acknowledgements}

Sandi Klav\v{z}ar acknowledges the financial support from the Slovenian Research Agency ARIS (research core funding P1-0297 and projects N1-0285, N1-0355).

%%%%%%%%%%%%%%%%%%%%%%%%%%%%
\section*{Declaration of interests}
%%%%%%%%%%%%%%%%%%%%%%%%%%%%
 
The authors declare that they have no conflict of interest. 

%%%%%%%%%%%%%%%%%%%%%%%%%%%%
\section*{Data availability}
%%%%%%%%%%%%%%%%%%%%%%%%%%%%
 
Our manuscript has no associated data.

%%%%%%%%%%%%%%%%%%%%%%%%%%%%%%%%%%%%%%%%%%%%%%%%%%%%%%

\end{document}